\documentclass{conm-p-l}
\usepackage{amssymb,graphics,apalike}
\newcommand{\y}{\\[5pt]}
\newcommand{\q}{\quad}
\newcommand\vs{\vskip5pt}\newcommand\bb{\bigbreak}
\newcommand\n{\noindent}\newcommand\lar{\leftarrow}
\newcommand\mb{\medbreak}\renewcommand\sb{\smallbreak}

\def\FIG{\vspace{-25pt}\centerline{\scalebox{.65}{\includegraphics{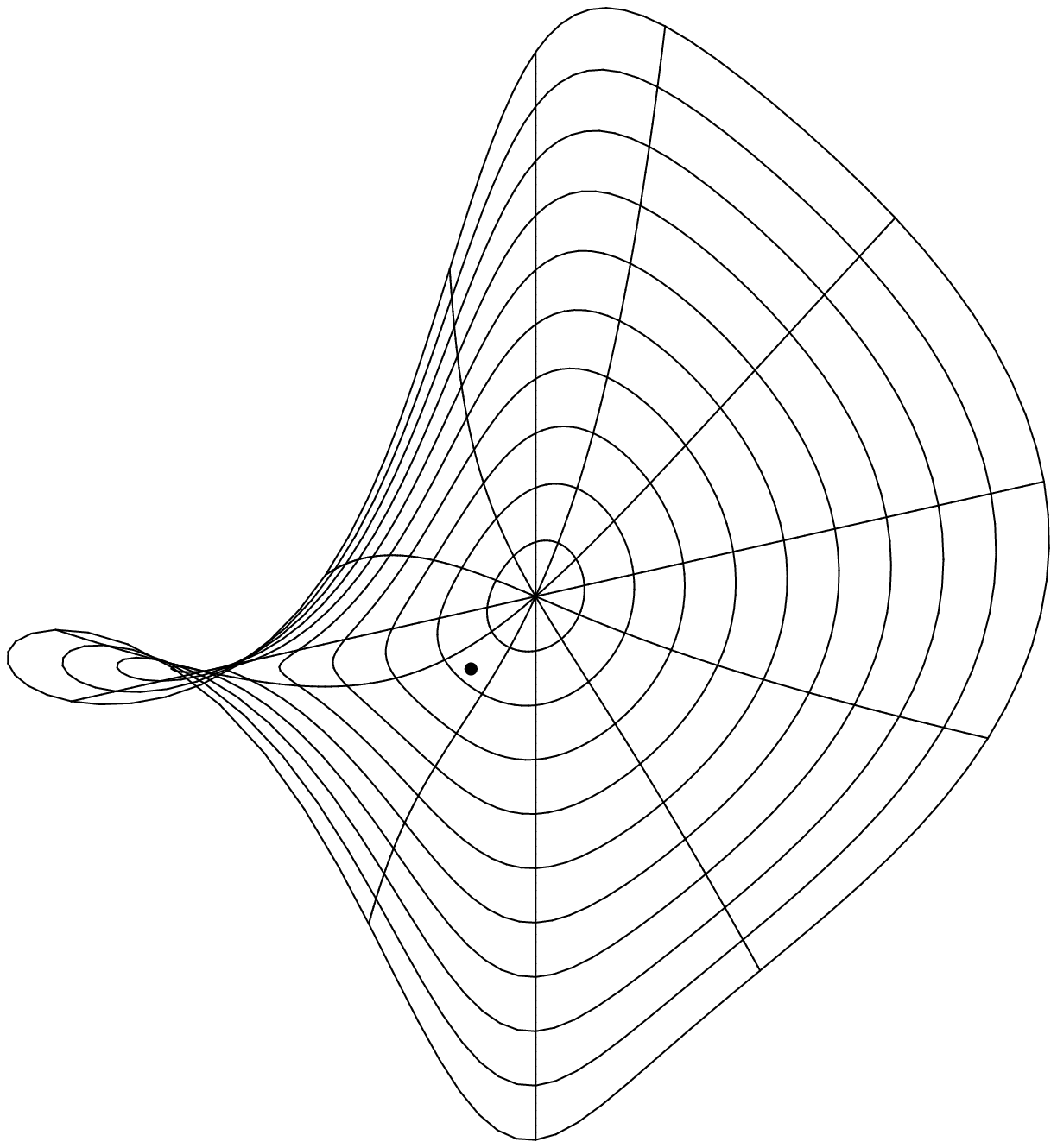}}}
\vspace{-15pt}\centerline{Figure 1}\vspace{25pt}}
\def\FIGS{\vfil\eject\phantom.\par
\centerline{\scalebox{.54}{\includegraphics{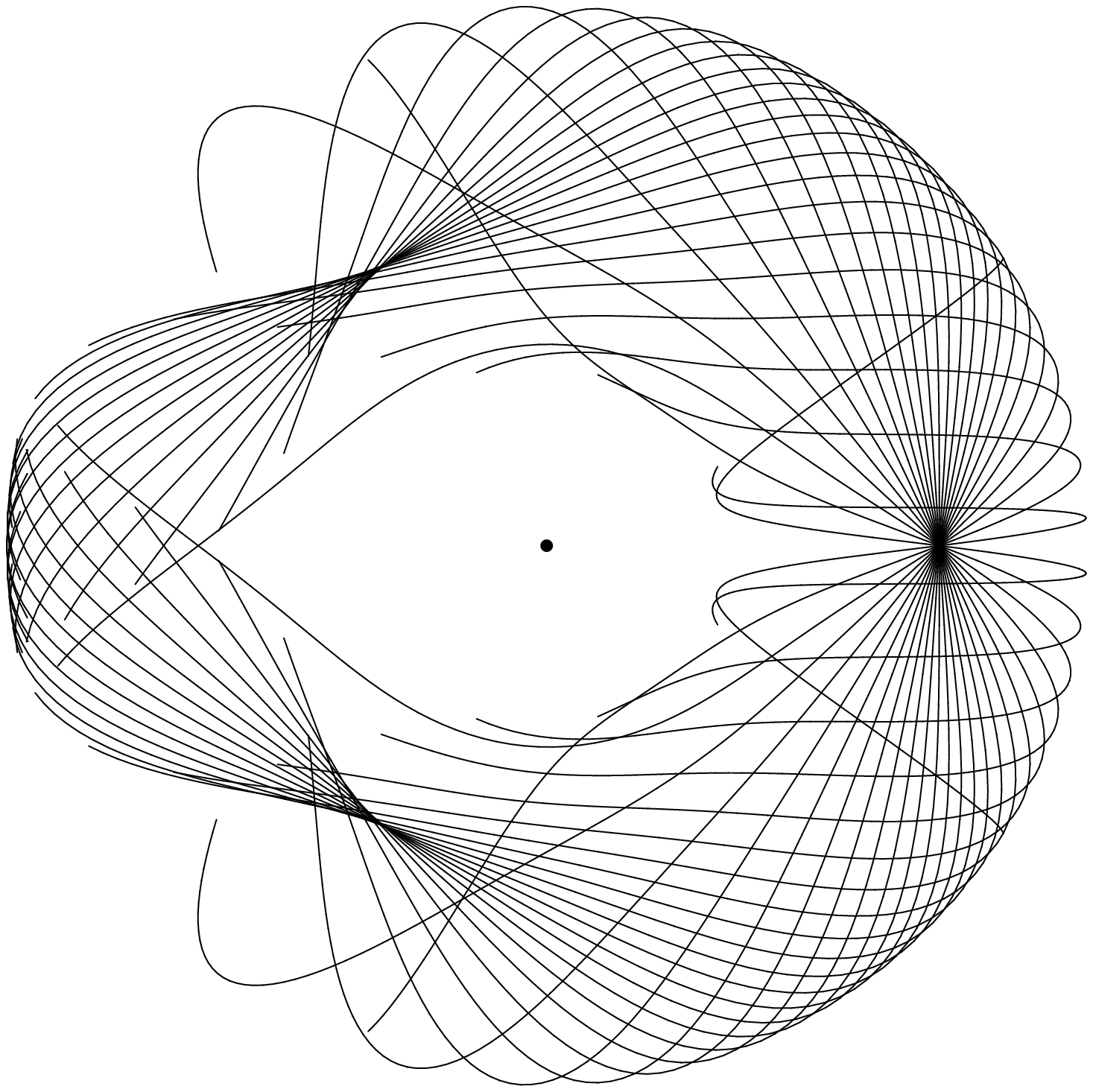}}}
\vspace{0pt}\centerline{Figure 2}\par\vspace{20pt}
\centerline{\hphantom{ooooo}\scalebox{.45}{\includegraphics{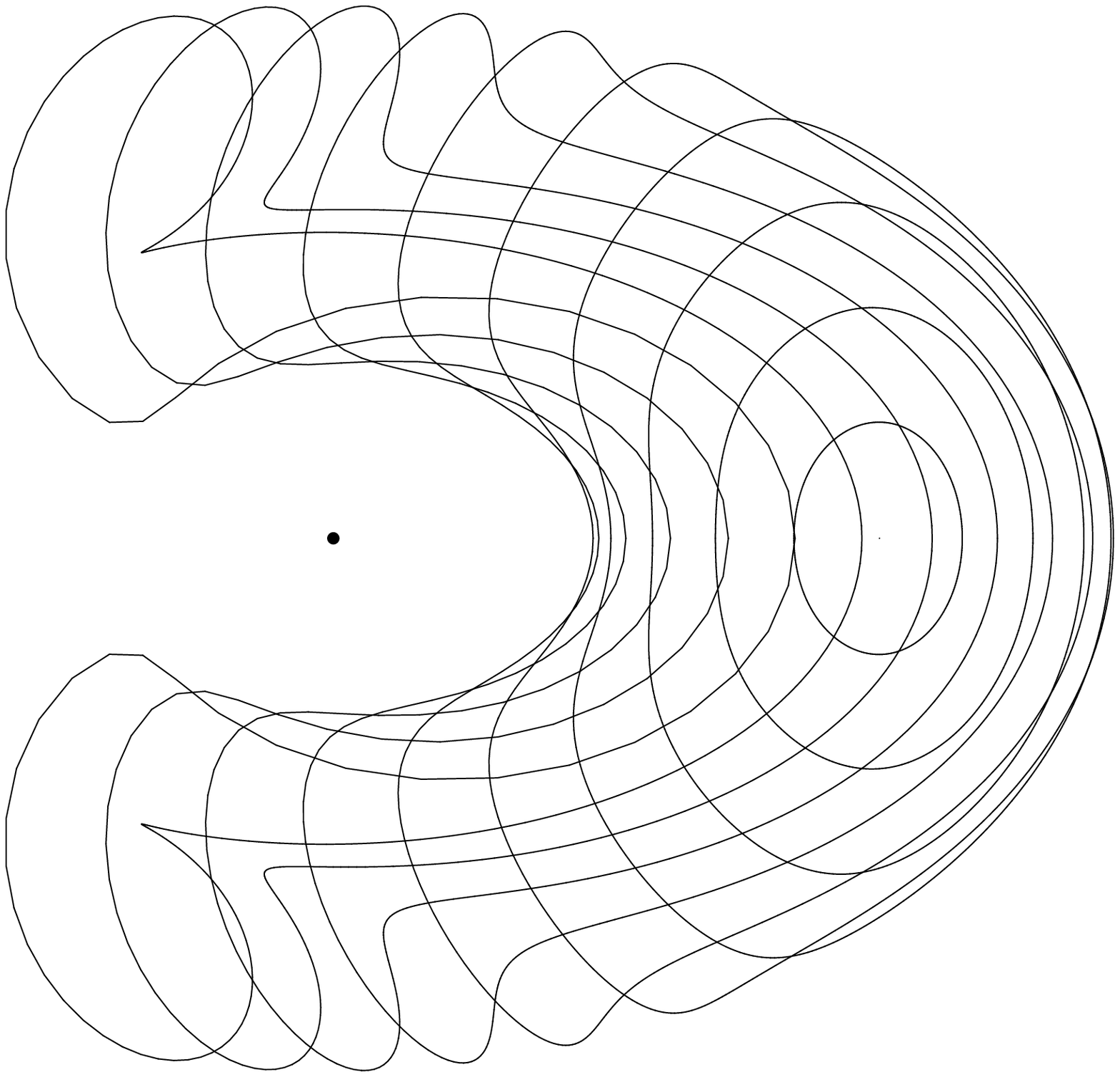}}}
\vspace{0pt}\centerline{Figure 3}\vfil\eject}
\newcommand{\ds}{\displaystyle}
\renewcommand\ge{\geqslant}
\newcommand{\g}{\mathfrak{g}}\newcommand{\h}{\mathfrak{h}}
\newcommand{\fm}{\mathfrak{m}}
\newcommand{\Ric}{\mathit{Ric}}\newcommand{\Wey}{\mathit{Wey}}
\newcommand{\Ga}{\Gamma}\newcommand{\Si}{\Sigma}
\newcommand{\Z}{\mathbb{Z}}
\newcommand{\C}{\mathbb{C}}\newcommand{\R}{\mathbb{R}}
\newcommand{\B}{\mathbb{B}}\newcommand{\CP}{\mathbb{CP}}

\renewcommand{\a}{\alpha}\renewcommand{\b}{\beta}
\newcommand{\na}{\nabla}\newcommand{\si}{\sigma}
\newcommand{\La}{\Lambda}\newcommand{\Spin}{\mathit{Spin}\kern1.5pt}
\newcommand\w{\omega}\newcommand\W{\Omega}\newcommand{\hW}{\widehat\W}
\newcommand{\we}{\wedge}
\newcommand{\op}{\oplus}\newcommand{\ot}{\otimes}
\newcommand{\rf}[1]{(\ref{#1})}\def\lie#1({\mathfrak{#1}(}
\newcommand{\ba}{\begin{array}}\newcommand{\ea}{\end{array}}
\def\be#1\ee{\begin{equation}#1\end{equation}}
\newcommand{\bt}{\begin{tabular}}\newcommand{\et}{\end{tabular}}
\newcommand{\al}{\langle}\newcommand{\ar}{\rangle}

\newcommand{\cR}{\mathcal{R}}\newcommand{\cW}{\mathcal{W}}
\def\stc{\stepcounter{equation}\ \textbf{\theequation}\ }
\def\pro#1#2\par{\bb\n\textbf{#1}\stc{\it#2}\par\bb}
\def\proc#1#2#3\par{\bb\n\textbf{#1}\stc\textrm{\cite{#2}} 
                    {\it#3}\par\bb}
\renewcommand\iff{\hbox{ $\Leftrightarrow$ }}
\renewcommand\t{\times}
\def\E{\raise1pt\hbox{$\bigwedge$}}
\newcommand{\ft}[2]{\hbox{$\textstyle\frac{#1}{#2}$}}
\newcommand{\fs}[2]{\hbox{\large$\frac{#1}{#2}$}}
\newcommand{\sss}[1]{\stepcounter{enumii}
                     \setcounter{equation}0\section{#1}}

\renewcommand\theequation{\arabic{enumii}.\arabic{equation}}
\newcommand\ex{\n\textit{Example}}
\begin{document}\parskip1pt\parindent15pt

\title{Almost Parallel Structures}

\author[S.\,M.~Salamon]{Simon Salamon}
\address{Department of Mathematics, Politecnico di Torino,
Corso Duca degli Abruzzi 24, 10129 Torino, Italy}
\email{salamon@calvino.polito.it}

\begin{abstract} A discussion of torsion of Riemannian $G$-structures leads 
to a survey of contributions of Alfred Gray and others on almost Hermitian
manifolds, $G_2$-manifolds, curvature identities, volume expansions, plotting
geodesics, and the geometry of nilmanifolds. The paper concludes with a new
example of a compact 8-manifold with a quaternionic 4-form that is closed but
not parallel.\end{abstract}

\maketitle

\section*{Introduction}

This article is closely based on the lecture entitled `A Mathematical Tribute'
that I gave at the conference in memory of Alfred Gray in Bilbao in September
2000. In preparing the lecture, I was well aware that a large number of his
theorems had directly influenced my own work. The aim therefore was to present
not just a personal selection of his results, but also an indication of
progress that has been made in many of the areas that he pioneered. In this
written version, I have given more emphasis to subsequent work, and have added
a final section that suits the new title. It is evident from a comparison with
the bibliography at the beginning of this volume that the topics below cover a
fraction of Gray's overall output. Nonetheless they demonstrate the importance
his work has had for the development of differential geometry.

I first met Alfred Gray at the IHES on Friday 13 April 1979, whilst visiting
Paris with fellow graduate student Martin Guest. I remember looking forward
with anticipation to a meeting on a combination of day and date that in the
past had seemed to bring me good fortune, and this occasion was to live up to
expectation. It led to my arrival in Maryland that summer, and regular contact
with Alfred that lasted until his untimely death in 1998. Although we never
formally collaborated, mathematical ideas were exchanged at dozens of memorable
occasions and locations.

Whilst finishing my thesis on quaternionic geometry, I had been working on
holonomy groups in general and studying properties of various hypothetical
geometries. At that time, there was a fair amount of pessimism about the
existence of the more exotic structures that Berger's list had
highlighted. Alekseevsky and Brown--Gray had supplied independent proofs that
any metric with holonomy $\Spin9$ is necessarily locally symmetric, and there
appeared to be great difficulty in constructing jets of metrics with holonomy
equal to $G_2$ or $\Spin7$. Equally disheartening was the belief that there
were no compact irreducible hyperk\"ahler manifolds in dimension greater than
4.

Against this background, Gray had been systematically examining $G$-structures
associated with potential holonomy groups, and paying close attention to
conditions that generalize the full holonomy reduction. I myself understood his
approach in terms explained in \S1, a discussion of which probably formed the
basis of our first meeting. Although some of the comments below now seem na\"\i
ve in the light of a greater understanding of exceptional geometry, everything
was less obvious twenty years ago.

\sss{Tensors and representations}

Let $M$ be an oriented Riemannian manifold of dimension $d$, with metric tensor
$g$ and Levi\,Civita connection $\na$. Let $\phi$ be a tensor with associated
stabilizer \[G=\{a\in SO(n):a\cdot\phi= \phi\},\] relative to a fixed
orthonormal frame at a fixed point $x\in M$. The restricted holonomy group of
$M$ is a subgroup of $G$ if and only if $\na\phi\equiv0$.

The space $\E^2T^*=\E^2T^*_xM$ of 2-forms at $x\in M$ is isomorphic to the Lie
algebra $\lie so(d)$ of skew-symmetric matrices. This contains the Lie algebra
$\g$ of $G$ and, relative to our fixed frame, there is a consistent splitting
\be\E^2T^*=\g\op \g^\perp.\label{gg}\ee
This can be used to provide a unified
description of the tensor $\na\phi$ that measures the failure of the holonomy
algebra to lie in $\g$.

There is no obstruction to the existence of a connection $\na'$ that preserves
the $G$-structure. If we choose one, then $\na_X-\na'_X$ is a tensor with
values in \rf{gg}, and $\na_X\phi=(\na_X-\na'_X)\phi$ belongs to $\g^\perp$ at
$x$. Thus,

\pro{Lemma} $\na\phi$ can be identified with an element of the space
$T^*\ot\g^\perp$.

\n It is easy to check that $\na\phi$ carries the same information as that
component of the torsion of $\na'$ that is independent of the choice of
$G$-connection. This component was described as the space of `intrinsic
torsion' in a similar context \cite{Bry}.\sb

For a given $G$, there will be a decomposition \be T^*\ot\g^\perp=
\bigoplus_{i=1}^N\cW_i,\label{dec}\ee in which each $\cW_i$ is an irreducible
$G$-module. Determining this and similar decompositions is a relatively
straightforward matter that nowadays can be performed quickly by computer.
Given his later admiration for computational methods, it is perhaps ironic that
Alfred carried out all the decompositions by hand. But in this way, he was able
to write down formulae for the various projections and tensors in meticulous
detail, and develop properties of the corresponding classes of manifolds that
others had overlooked.\sb

\ex s (i) $d=2n$, $\phi=J$ is an orthogonal almost complex structure, and
$G=U(n)$. In this case, $N=4$ provided $n\ge3$. The same theory results if
$\phi=\w$ is a non-degenerate 2-form compatible with $g$. Indeed, 3 of the 4
real components of $T^*\ot\g^\perp$ can be identified with those of the
exterior derivative \[d\w\in \E^3T^*=[\![\La^{3,0}]\!]\op[\![\La^{2,1}_0]\!]\op
T^*,\] notation as in \cite{Shol}. The fourth component forms part of the
Nijenhuis tensor of $J$. Given an even-dimensional Riemannian manifold $(M,g)$,
one may investigate the existence of compatible (i) complex structures, and
(ii) symplectic structures, and there are constraints arising from the
curvature tensor of a rather different nature in each of the two cases.\sb

\n(ii) $d=7$, $\phi$ is a `positive generic 3-form', meaning that there exists
an orthonormal basis of 1-forms at each point relative to which
\be\phi=e^{125}- e^{345}+e^{136}-e^{426}+e^{147}-e^{237}+e^{567},\label{phi}\ee
and $G$ is the exceptional Lie group $G_2$. As was well understood from the
study of vector cross products \cite{GX}, the 3-form $\phi$ can also be defined
in terms of Cayley multiplication on $\hbox{Im}\,\mathbb{O}=\R^7$. Once again,
$N=4$. The four components of $T^*\ot\g^\perp$ can be recovered from those of
$d\phi$ and $d*\phi$ as explained in \S3, and the metric has holonomy contained
in $G_2$ if and only if both these forms vanish. The larger class of
`calibrated' $G_2$-manifolds consists of ones for which $d\phi=0$; by contrast
$d*\phi$ is more akin to the Nijenhuis tensor in case (i).\sb

\n(iii) $d=4n$, $\phi$ is a `quaternionic' 4-form, meaning that \[\textstyle
\phi=\sum\limits_1^3\w^i\!\we\w^i,\] where $\w_1,\w_2,\w_3$ is a triple of
2-forms corresponding to almost complex structures $I_1,I_2,I_3$ modelled on
the Lie algebra $\lie su(2)$. The stabilizer of $\phi$ is $Sp(n)Sp(1)$, and
such structures were considered in \cite{Gnote}. Whilst he did not study the
full decomposition of $T^*\ot\g^\perp$, it is known that $N=6$ provided that
$n\ge3$ \cite{Shol}. The vanishing of exactly 3 of these 6 components
characterizes the class of quaternionic manifolds which possess a complex
twistor space.\mb

We shall have something to say about each of these structures in
\S2,\,\S3,\,\S7 respectively. Other examples studied by others include $\Spin7$
structures (for which $N=2$), almost product structures ($N=6$), and almost
contact structures. See also the analysis by Friedrich of the $\Spin9$ case in
this volume. We shall confine ourselves to two additional remarks before
leaving the general set-up.

Firstly, the decomposition \rf{gg} has been exploited to define a type of
Dolbeault complex for $G$-structures, important for the study of instantons
associated to special geometries \cite{RC}. Let $N(G)$ denote the normalizer of
$G$ in $SO(d)$; if $N(G)$ is strictly larger than $G$, it suffices to start
with an $N(G)$-structure in order to define \rf{gg}. On can then set
$A^0=\E^0T^*$, $A^1=T^*$, $A^2=\g^\perp$, and more generally \[A^k=
(\g\we\E^{k-2}T^*)^\perp,\q k\ge 3.\] Let $d_k$ denote the composition of
exterior differentiation acting on $A^k$ with projection $\E^{k+1}T^*\to
A^{k+1}$. An `extendability' condition guarantees that \[A^0\to A^1\to A^2\to
A^3\to\cdots\] is a complex, and can be expressed in terms of the vanishing of
certain components of $\cW$. The case $G=G_2$ was developed by \cite{Ug}, based
on a thesis co-examined by Gray and myself.

Secondly, relative to an orthonormal frame compatible with the holonomy group
$G$, the Riemann curvature tensor has components $R_{ijkl}$ whose symmetries
determine an element in the space \be \cR=\ker\Big(S^2\g\hookrightarrow
S^2(\E^2T^*)\to \E^4T^*\Big),\label{bia}\ee in which the second mapping
corresponds to the first Bianchi identity. Thus, any $G$-invariant element $B$
(such as the Killing form) in $S^2\g$ will have non-zero image in $\E^4T^*$
unless it satisfies the same symmetries as a Riemann curvature tensor. It
follows from the theory of E.~Cartan that any $G$-structure with $G\subset
SO(n)$ possesses a canonical 4-form, provided the structure possesses metrics
that are not locally symmetric. Example (iii) above is actually defined by such
a 4-form, as is $\Spin7$ geometry, though invariant 4-forms are also present in
(i) (namely $\w\we\w$), and in (ii) (namely $*\phi$).

\sss{The sixteen classes of almost Hermitian manifolds}

Let $M$ be an almost Hermitian manifold of dimension $d=2n$. It is therefore
equipped with a positive-definite metric tensor $g$, an orthogonal almost
complex structure $J$, and associated 2-form $\w$. The subgroup of $SO(2n)$
determined by these structures is isomorphic to the unitary group $U(n)$.

The structure is K\"ahler if and only if $\na J\equiv0$, and the decomposition
\rf{dec} in this case is given by

\proc{Proposition}{GH80} For $n\ge3$, the tensor $\na J$ belongs to the direct
sum \be\cW_1\op\cW_2\op\cW_3\op\cW_4\label{W4}\ee of four irreducible
$U(n)$-modules of respective dimensions
\[\ft13n(n-1)(n-2),\q\ft23n(n-1)(n+1),\q n(n+1)(n-2),\q 2n.\]\vs

The list of dimensions reflects the following facts:\sb

\n(i) When $n=2$, $\cW_1$ and $\cW_3$ reduce to zero, leaving the sum
$\cW_2\op\cW_4$ of two spaces each of dimension 4. The 2-form is closed if and
only if $\na J\in\cW_2$, and $J$ is integrable if and only if $\na J\in\cW_4$.
Thus, in 4 dimensions the symplectic and complex conditions are exactly
complementary. A fuller description of this case is included in the paper by
Donaldson in these Proceedings.\sb

\n(ii) The space $\cW_4$ is isomorphic to the cotangent space $T^*$. Such a
summand is always present in \rf{dec}, as it reflects the change in $\na\phi$
that occurs when the metric is altered by a conformal factor. For example, if
$M$ is locally conformally K\"ahler then necessarily $\na J\in\cW_4$.
Conversely, if $\na J\in \cW_4$ then $\na J$ can be identifed with a 1-form
$\theta$ that satisfies $d\w=\theta\we\w$. If $n\ge3$ then $d\theta=0$ and $M$
is locally conformally K\"ahler, but the 4-dimensional case is more
complicated. Of particular interest are those Hermitian 4-manifolds for which
$\na\theta=0$.\sb

\n(iii) The space $\cW_1$ underlies the third exterior power of the space of
$(1,0)$-forms. When $n=3$, $\cW_1$ has dimension 2, and a non-zero element of
it determines a reduction from $U(3)$ to $SU(3)$. This leads to some special
features of so-called nearly-K\"ahler manifolds of dimension 6 (see below).\sb

Further to (i), the two halves $\cW_1\op\cW_2$ and $\cW_3\op\cW_4$ of equal
dimension $n^2(n-1)$ lead to a sort of duality. They are characterized by the
conditions \[\ba{l} \na J\in\cW_1\op\cW_2 \iff (d\w)^{1,2}=0,\y\na J\in\cW_3
\op\cW_4 \iff M\hbox{ is Hermitian}.\ea\] The significance of the first
condition was already recognized by Gray in his thesis, and the manifolds
satisfying $(d\w)^{1,2}=0$ are often called \textit{quasi-K\"ahler}.

\proc{Proposition}{G65} If $M$ is quasi-K\"ahler, and $M'\subset M$ is a
$J$-holomorphic submanifold, then $M'$ is minimal.

\n This generalizes the well-known fact that a complex submanifold of a
K\"ahler manifold is minimal. A similar result, but expressed in terms of
harmonic mappings was presented by \cite{L}.\sb

The following case has rightly merited special attention.

\pro{Definition} $M$ is \textit{nearly-K\"ahler} if $\na J\in\cW_1$,
equivalently $(\na_XJ)X=0$ for all $X$. $M$ is said to be \textit{strictly}
nearly-K\"ahler if in addition $\na J\ne0$.

\n The basic model of a strictly nearly-K\"ahler manifold is the sphere
$S^6=G_2/SU(3)$ endowed with a non-integrable almost-Hermitian structure
induced from the cross product on $\hbox{Im}\,\mathbb{O}=\R^7$. In fact, $S^6$
belongs to a large class of homogeneous examples that we describe next.

Let $G$ be a real semisimple Lie group. A homogeneous manifold $M=G/H$ is
called a \textit{3-symmetric space} if $H$ is the fixed point set of an
automorphism $\theta$ of $G$ with $\theta^3=1$. The complexified Lie algebra of
$G$ has the form \be\g=\h\op\fm,\q \fm_c=\fm^{1,0}\op\fm^{0,1},\label{fm}\ee
and the resulting 3 summands of $\g_c$ are the eigenspaces of
$\theta$. Defining \be J=\fs1{\sqrt3}(2\theta+1)\label{J}\ee determines a
canonical almost complex structure on the tangent space $T_xM\cong\fm$. This is
non-integrable unless $M$ is Hermitian symmetric, a possibility we exclude
below.

\proc{Theorem}{WG} Any 3-symmetric space has a nearly-K\"ahler metric
compatible with \rf{J}.\sb

\n\textit{Classification.} The collection of 3-symmetric spaces can be divided
into the following categories, properties of which appear in \cite{G72}.

\n(i) generalizations of $S^6$ with irreducible isotropy, to cite one
example $E_8/SU(9)$. Such spaces are in a sense rarer than their
complex homogeneous counterparts.

\n(ii) $G\t G$, where $G$ is a compact Lie group. Just as such a group can be
viewed as an ordinary symmetric space by considering the space of cosets $(G\t
G)/G$, so $G\t G$ becomes a 3-symetric space when viewed as $(G\t G\t G)/G$.

\n(iii) the total space of a fibration $\pi:Z\to M$ where $M$ is a symmetric
space defined by an inner involution. In this case, $\fm^{1,0}$ is a reducible
representation of $H$, reflecting the horizontal/vertical decomposition of
$TZ$.\sb

The `twistor spaces' $Z$ in (iii) have been extensively used in the study of
minimal surfaces and harmonic maps, in essence because of results
\cite{BryL,S1164} that extend Proposition~2.3. For example,

\pro{Proposition} If $f\colon\Si\to Z$ is a $J$-holomorphic curve in a
3-symmetric space then $\pi\circ f$ is a harmonic map into $M$.\sb

\ex. A `classical' instance of this construction is provided by the diagram
\be\ba{rcc}&&\ds\frac{U(p+q+1)}{U(p)\t U(q)\t U(1)}=Z\\[8pt]&
\nearrow&\;\Big \downarrow\>\pi\\[8pt]\Si&\to &\ds\frac{U(n+1)}{U(n)\t
U(1)}=\CP^n\ea\label{Z}\ee

\n in which $p+q=n$. Any twistor space $Z$ has an integrable complex structure
$J_1$ in addition to $J_2$, and the two structures coincide on the horizontal
subspace of each tangent space $T_xZ$. The most important class of maps into
$Z$ are those that are horizontal and (unambibuously) holomorphic. If
$\Si\cong\CP^1$ has genus 0, any harmonic map into $\CP^n$ is the projection of
a horizontal holomorphic one in a flag manifold such as $Z$ above \cite{EW}.\sb

Related methors have been successfully used to classify orthogonal complex
structures and to extend the description of harmonic maps to symmetric spaces
\cite{BR,BuG}.

\sss{Weak and exceptional holonomy}

The importance of nearly-K\"ahler manifolds is highlighted by

\proc{Theorem}{GnK} A 6-dimensional strictly nearly-K\"ahler manifold
is Einstein.

\ex. It is well known that there are two Einstein metrics on $\CP^3$ for which
the Penrose twistor fibration $\pi:\CP^3\to S^4$ is a Riemannian
submersion. The non-standard one is nearly-K\"ahler relative the non-integrable
almost complex structure $J_2$ associated to $\pi$.\sb

Further curvature identities amount to the statement that the curvature of a
strictly nearly-K\"ahler 6-manifold has the form \[R=sR_{S^6}+R_{CY},\] where
$s$ is the scalar curvature, $R_{S^6}$ is the (suitably normalized) constant
curvature tensor of $S^6$, and $R_{CY}$ stands for the curvature tensor of a
6-manifold with a Calabi-Yau metric, so holonomy equal to $SU(3)$. The only
known compact 6-dimensional strictly nearly-K\"ahler manifolds are \[S^6,\q
S^3\t S^3,\q \CP^3,\q \mathbb{F}^3.\] The second arises from (ii) in \S2 above
with $G=SU(2)$, and the last one denotes the flag manifold $Z$ in \rf{Z} with
$p=q=1$. Any 6-dimensional strictly nearly-K\"ahler manifold has a canonical
connection with holonomy $SU(3)$, and this reduces to $U(2)$ if and only if $M$
fibres over a self-dual Einstein 4-manifold \cite{BM}. This construction is
important because it can be used to show that there is an abundance of
locally-defined strictly nearly-K\"ahler metrics.

Nearly-K\"ahler manifolds play a central role in the theory of weak holonomy
developed in \cite{G71}. A full holonomy reduction to $G$ requires a tensor
with stabilizer $G$ to remain constant under parallel translation, whereas $G$
is a weak holonomy group if a certain family of subspaces is mapped into itself
by parallel translation. An obvious candidate would be the family of
quaternionic lines (which are real 4-dimensional subspaces) for $G=Sp(n)Sp(1)$,
but it turns out that in this case a weak holonomy reduction is equivalent to a
full holonomy reduction. Instead, the theory is most fruitful in 6 and 7
dimensions.

Manifolds with weak holonomy $G_2$ can be characterized algebraically, as in
the nearly-K\"ahler case which corresponds to weak holonomy $U(3)$. If a
manifold $M$ has a 3-form of type \rf{phi}, then there are $G_2$-equivariant
decompositions \be\ba{l} \E^2T^*=\g\op\g^\perp \>\cong\> \g_2\op T^*\y
\E^3T^*\>\cong\>\R\op T^*\op S^2_0T^*.\ea\label{ext}\ee The presence of a
summand $T^*$ in $\E^2T^*$ corresponds to the cross product operation, and the
presence of $T^*$ in $\E^3T^*$ reflects the failure of this product to satisfy
the Jacobi identity. The decomposition \rf{dec} becomes \be \ba{rcl}
T^*\ot\g_2^\perp\>\cong\> T^*\ot T^* &\cong& S^2T^*\op\E^2T^*\y &\cong&\R\op
S^2_0T^*\op T^*\op\g_2,\ea\label{4}\ee and a comparison of \rf{ext} and \rf{4}
reveals

\proc{Corollary}{FG82} $M$ has holonomy in $G_2$ iff $d*\phi=0$ and
$d\phi=0$.\sb

The manifold has weak holonomy $G_2$ if and only if $\na\phi$ lies in the
1-dimensional component of \rf{4}. This means that \be
d\phi=c*\phi\label{3*4}\ee for some $c$ (a universal constant times the scalar
curvature that we assume is non-zero), and consequently $d*\phi=0$. In this
case, associative subspaces are ones invariant by the cross product, and are
preserved by parallel transport.  \sb

The local existence of metrics with holonomy equal to $G_2$ was proved
by \cite{Bry}. The first space found to exhibit an explicit metric with
exceptional holonomy was the cone $Z=Y\t\R^+$, in which
\[Y=\frac{SO(5)}{SO(3)}\] is an isotropy irreducible space arising
from the representation $S^6\C^2$ of $SO(3)\!=\!SU(2)/\Z_2$, with a history of
providing counterexamples \cite{Ber}. The fact that there is a 1-dimensional
space of invariant 3-forms in $\E^3(S^6\C^2)$ led the author to observe that
$Y$ has a metric with weak holonomy $G_2$ satisfying \rf{3*4}. The realization
that $Z$ has a closed 4-form defining a metric with holonomy equal to $\Spin7$
was a step away. There is no doubt that the discovery of explicit metrics with
exceptional holonomy would have taken longer without Gray's theory of weak
holonomy. The subsequent theory of compact manifolds with exceptional holonomy
is the subject of work culminating in the book \cite{Joy}.

With the hindsight provided by the first examples of metrics with exceptional
holonomy, it became clear that the most effective way of understanding the
relationship between weak and exceptional holonomy is via the theory of Killing
spinors, developed by \cite{Fetc}. Indeed,

\proc{Theorem}{Bar} (i) A Riemannian metric on $X^6$ is strictly
nearly-K\"ahler iff the corresponding conical metric on $X\t\R^+$ is Ricci-flat
with holonomy contained in $G_2$.\y (ii) A metric on $Y^7$ has weak holonomy
$G_2$ iff the conical metric on $Y\t\R^+$ is Ricci-flat metric with holonomy
contained in $\Spin7$.

\n Possible holonomy groups of $Y\t\R^+$ in (ii) are $\Spin7$, $SU(4)$ and
$Sp(2)$, and these correspond to 7-manifolds with a space of Killing spinors of
dimension 1,2,3 respectively. This is fully discussed in \cite{FKMS}, and
examples appear in \cite{GS}.\mb

\sss{Curvature and volume}

Gray's published papers incorporate extensive studies of the Riemann curvature
tensor. As I realized during the Bilbao conference (responding to a query from
Peter Gilkey), given any natural condition on the curvature, the chances are
that Alfred had already given it a name. His classification of manifolds in
terms of the covariant derivative of the Ricci tensor \cite{GEin} is
particularly appealing, and was the subject of talks at the
conference. However, in this section, we have chosen to comment on the
curvature of Hermitian manifolds, and unrelated volume-type expansions.\sb

A striking result concerning the curvature of K\"ahler manifolds is

\proc{Theorem}{G77} A compact K\"ahler manifold with nonnegative sectional
curvature and constant scalar curvature is locally symmetric.

\n This result is a prototype for many subsequent theorems in the literature,
characterizing the curvature of Hermitian symmetric spaces and complex
projective space.

The curvature tensor $K$ of a K\"ahler manifold of real dimension $2n$
satisfies the equation \be K(W,X,Y,Z)=K(W,X,JY,JZ)\label{Kah}\ee that reflects
the fact that $K$ may be regarded as a 2-form with values in the holonomy
algebra $\lie u(n)$. The curvature tensor $R$ of an arbitary almost Hermitian
manifold can therefore be decomposed as \[R=K+K^\perp,\] where $K$ represents
the component of $R$ in the subspace of tensors satisfying \rf{Kah}. The next
result (see for example \cite{FFS}) is a starting point for the analysis of
metrics with constant holomorphic sectional curvature.

\pro{Proposition} The tensor $K^\perp$ is a linear contraction of
$\na\na J$ (where $\na$ denotes the Levi\,Civita connection) and has
zero holomorphic sectional curvature.\sb

There are further decompositions of $K$ and $K^\perp$ under the action of
$GL(n,\C)$, and decompositions under $U(n)$ are given in \cite{TV}. A
significant fact concerning curvature of a Hermitian manifold is encapsulated
in

\proc{Proposition}{G76} If $M$ is Hermitian then \[\ba{l}R(W,X,Y,Z)+
R(JW,JX,JY,JZ)\\\hspace{40pt}=R(JW,JX,Y,Z)+R(JW,X,JY,Z)+R(JW,X,Y,JZ)\\
\hspace{70pt}+R(W,JX,JY,Z)+R(W,JX,Y,JZ)+R(W,X,JY,JZ).\ea\]

\n A complex structure that is orthogonal relative to a metric $g$ is also
orthogonal relative to any conformally related metric. It follows that the
above equation constrains only the Weyl tensor $\Wey$ of $M^{2n}$. In terms of
type decomposition it amounts to saying that $R$ has no component in the real
subspace underlying $S^2(\La^{2,0})\ominus\La^{4,0}$, and imposes $k=\frac16
n^2(n^2-1)$ equations on $\Wey$, which itself has dimension asymptotic to $8k$
as $n\to\infty$.

The Weyl tensor of an oriented Riemannian 4-manifold $M$ decomposes as \be
\Wey=\Wey_++\Wey_-,\q\Wey_\pm\in\E^2_\pm T^*M.\label{Wpm}\ee If $\Wey_+\ne0$
then $M$ has at most two distinct pairs $\pm I_1,\,\pm I_2$ of orthogonal
complex structures compatible with the metric and orientation. In higher
dimensions, less is known about metrics with a `large' number of isolated
orthogonal complex structures, though work of \cite{Pon,Apo} suggests that the
theory will have both a local and global flavour.\mb

The curvature of a Riemannian manifold provides a quantifiable means of
comparing the volume of submanifolds with those of Euclidean space, and
established theorems relate the volume of balls and spheres in Riemannian
manifolds with those of spaces of constant curvature \cite{Gun,BC}. Gray
extended this theory by developing asymptotic expansions to study the volume of
balls and subsequently tubes around curves and submanifolds.

Let $V(r)$ denote the volume of a ball $\B(r)$ of radius $r$ centred at the
origin in $\R^d$. It is convenient to set $d=2n$, whether or not $d$ is even.
Then \[V(r)=V(1)r^d=\frac{(\pi r^2)^n}{n!},\] where $n!$ stands for $\Ga(n+1)$
if $n$ is half-integral. The theory of volume in higher dimensions incorporates
a number of counter-intuitive features that lend themselves naturally to
investigation in curved space:\sb

\n(i) The table below shows that $V(1)$ reaches a peak for $d=5$ and then
decreases to 0. \[\ba{c||c|c|c|c|c|c|c|c|c|c|c|c|} d & 1&
2&3&4&5&6&7&8&9&10&15&20\\\hline V(1)&
2&3.14&4.19&4.93&5.26&5.17&4.72&4.06&3.30&2.55&0.38&0.03\ea\] To express this
in another way, let $r_d$ be the radius for which $V(r_d)=1$. For instance,
$r_{1000}\sim7.68$, and Stirling's formula implies that \[ r_d\sim\sqrt{\frac
d{2\pi e}},\q d\to\infty.\]

\n(ii) A study of the distribution of volume within balls implies that
(for example) the slice \[\B(r_d)\cap(\R^{d-1}\t[-0.4,0.4])\] has
volume at least $0.8$ irrespective of $d$ (as the author discovered
whilst browsing the conference bookstall). This is a simple result
concerning a `tube' around a planar hypersurface of $\B(r_d)$.\sb

Now suppose that $M$ is a Riemannian manifold. Fix a point $x\in M$, and let
$V(r)$ now denote the volume of a ball $\B(r)$ of radius $r$ centered at
$x$. We assume that $r$ is less than the injectivity radius of $M$, so that
$\B(r)$ is formed of geodesics of length $r$ emanating from $x$.  Early results
on the volume of $\B(r)$ include those of \cite{Hotel,Weyl}. There is an
expansion \[V(r)=\frac{(\pi r^2)^n}{n!}\left(1+c_2r^2+c_4r^4+c_6r^6+\cdots
\right),\] in which the coefficient $c_2$ equals $-\frac13s/(n+1)$ where $s$ is
the scalar curvature ($n=\frac12\dim M$). Gray took up the challenge of
determining further coefficients. Let $R$ denote the full curvature tensor
(with components $R_{ijkl}$ relative to an orthonormal frame), and $\Ric$ the
Ricci tensor (components $R_{jl}=\sum_i R^i_{jil}$).

\proc{Theorem}{G73} $\ds c_4=\frac{8\|\Ric\|^2-3\|R\|^2+5s^2-18\Delta s}
{1440(n+1)(n+2)}$.\sb

A natural conjecture is that the vanishing of all $c_{2k}$ implies that $M$ is
flat, but a proof of this became elusive by an array of examples of
$2n$-manifolds for which \be V(r)=\frac{(\pi r^2)^n}{n!}(1+O(r^{2k})). 
\label{ce}\ee Results from \cite{GV} showed that in
dimension $4$, there exist metrics with $s=0$ and $c_4=0$ so that \rf{ce} holds
with $k=3$. The same paper has an example with $n=367$ and $k=4$. These results
were drammatically extended by the use of an additive functor to prove

\proc{Theorem}{Kow} There exists a product of homogeneous spaces such that
\rf{ce} holds with $k=8$.\sb
 
Similar techniques can be applied to other asymptotic expansions involving
curvature. An analogous conjecture concerning the class of so-called harmonic
manifolds motivated the paper \cite{CGW}, but was settled in the negative by
the discovery of the so-called Damek-Ricci spaces. On a different topic, the
article by Pinsky in the Proceedings refers to work on `mean exit times' for
Brownian motion on Riemannian manifolds \cite{GPin}.\sb

The volume of tubes is the subject of many of Alfred Gray's papers that
provided the basis of his first book \cite{Gtubes}. It investigates a theory
that has developed from work of H.~Weyl, whose starting point is the fact that
the volume of a tube of radius $r$ about a space curve of arclength $\ell$
remarkably depends only on $r$ and $\ell$.

Corresponding invariance properties in the context of K\"ahler manifolds and
characteristic classes lead to a number of elegant formulae such as

\proc{Theorem}{G85} A tube of (sufficiently small) radius $r$
surrounding a hypersurface of degree $k$ in $\CP^n$ has volume
$\ds\frac{\pi^n}{n!}\Big(1-(1-k\sin^2\!r)^n\Big)$.\vs

\sss{Plotting geodesics on surfaces}

We have been talking about geodesic balls on manifolds. Such objects can be
graphically illustrated in the 2-dimensional case. The first edition of Gray's
book on curves and surfaces contained some of the first programs for
constructing geodesics on surfaces.

\FIG

By way of a diversion at the conference, I displayed a geodesic clockface on
the paraboloid $z=xy$, plotted using the \textsc{Mathematica} program
\texttt{solvegeoeqs} and its cousins \cite{Gbook}. Another is shown in
Figure~1, in which the saddle point is indicated by a dot and one can easily
spot two straight line rulings. Analogous curves on a torus are displayed in
Figures~2 and 3, and are not entirely inconsistent with the title of this
article. The aim of this section is to present a self-explanatory account of
how the plotting is quickly achieved, to serve as an introduction to the more
extensive programs in \cite{Gbook}. This powerful construction kit can be
downloaded from the website \texttt{math.cl.uh.edu/~gray} maintained by
M.~Mezzino.

Figure~2 consists of geodesics emanating from a point $P$ on a torus of
revolution with conjugate points visible, and Figure~3 displays corresponding
`circles' centered at $P$ up to the point where they become singular. They were
plotted with the program below, a condensed version of Gray's that can be
copied by hand onto a keyboard without too much effort. It mixes local and
global variables in a rather amateur way, but is designed to make the various
steps transparent. The instructions can easily be modified so as to plot
individual geodesics and related objects on an arbitrary surface.

Let $\mathbf{x}(u,v)$ be a parameterized surface. The aim is to construct
geodesics emanating from a fixed point $\mathbf{x}(a,b)$, and to this end we
divide the task into three blocks of code that correspond to notebook
cells. The first consists of standard definitions that are independent of the
choice of the function $\mathbf x$. The latter is entered at the second stage,
and plotting takes place after loading the required parameters.\mb

\n(i) First, one defines the coefficients of the first fundamental form and the
Christoffel symbols. These basic differential geometric formulae translate
readily into computer code, in readiness for all manner of applications.\sb
\begin{verbatim}
    Unprotect[E];
    E:= D[x[u,v],u].D[x[u,v],u]
    F:= D[x[u,v],u].D[x[u,v],v]
    G:= D[x[u,v],v].D[x[u,v],v]
    ga[1,1,1]:= D[E,v]F-2D[F,u]F+D[E,u]G
    ga[2,2,2]:= D[G,v]E-2D[F,v]F+D[G,u]F
    ga[2,1,1]:=-D[E,v]E+2D[F,u]E-D[E,u]F
    ga[1,2,2]:=-D[G,v]F+2D[F,v]G-D[G,u]G
    ga[1,1,2]:= D[E,v]G-D[G,u]F
    ga[2,1,2]:=-D[E,v]F+D[G,u]E
    Ga[i_,j_,k_]:= Simplify[ga[i,j,k]/(E G-F^2)/2]
\end{verbatim}\vs

\n(ii) It is convenient to enter the required parameterization at this
point. The geodesic equations are then entered, and are then in a form ready to
use for the surface in question. Definitions are included of appropriate
initial conditions and the command for solving the differential equations.\sb
\begin{verbatim}
    x[u_,v_]:= {(2+Cos[v])Cos[u],(2+Cos[v])Sin[u],Sin[v]}
    su:= {u->u[s],v->v[s],p->u'[s],q->v'[s]}
    e[j_]:=e[j]= Ga[j,1,1]p^2+2Ga[j,1,2]pq+Ga[j,2,2]q^2 /.su
    eqic:= {u''[s]+e[1]==0,v''[s]+e[2]==0,
            u[0]==a,v[0]==b,u'[0]==Cos[th],v'[0]==Sin[th]}
    so:= NDSolve[eqic,{u,v},{s,0,r}]
    sd[m_]:= Flatten[Table[so,{th,0,2Pi,2Pi/m}],1]
\end{verbatim}\vs

\n(iii) Points along geodesics equidistant from $(a,b)$ are joined and
plotted, together with the geodesics themselves. It suffices to re-enter the
lines below to re-draw the plot with different parameters.

\FIGS

\begin{verbatim}
    a:=0; b:=Pi/2; r:=4; m1:=100; m2:=50; k:=.4
    xs:= x@@Sequence[{u[t],v[t]}]
    g0:= Graphics3D[Point[{0,0,0}]]
    g1:= Graphics3D[Table[Line[xs/.sd[m1]],{t,0,r,k}]]
    g2:= ParametricPlot3D[Evaluate[xs/.sd[m2]],{t,0,r},
                          DisplayFunction->Identity]
    Show[g0,g1,g2,Boxed->False,BoxRatios->{1,1,1},ViewPoint->{0,0,4}]
\end{verbatim}\vs     

\n The graphics object \texttt{g0} enables the origin to be added to the
diagram as a point of reference. Other commands in (iii) that enable the
various definitions to be combined are of course fully explained in the
Mathematica Book \cite{Wol}.

\sss{Invariant geometry on nilmanifolds}

Any compact simple Lie group $G$ of even dimension $d=2n$ admits left-invariant
complex structures, that can easily be described in terms of the root
decomposition of the Lie algebra $\g$. On the other hand, $G$ admits no
symplectic structure as $b_2=0$. By contrast, a given nilpotent Lie group $N$
may or may not admit left-invariant complex and/or symplectic structures.\sb

\ex. There are three simply-connected nilpotent Lie groups in dimension $d=4$,
namely $N_1,N_2,N_3$ in which the Lie-algebra of $N_k$ is $k$-step (so $N_1$ is
abelian). All three admit left-invariant symplectic structures, but it is easy
to see that only the first two admit invariant complex structures.\sb

If $N$ is a nilpotent Lie groups with rational structure constants, it
possesses a discrete subgroup $\Gamma$ for which $M=N/\Gamma$ is compact
\cite{Mal}. Such a compact nilmanifold $N/\Ga$ can only admit a K\"ahler metric
if $N$ is abelian, in which case the quotient is a torus (see for example
\cite{BG}). A celebrated theorem \cite{Nom} implies that the space of
left-invariant forms provides a minimal model for the nilmanifold's de\,Rham
cohomology, which is therefore readily computed. Gray pioneered the use of
Massey products to detect the non-existence of a K\"ahler metric, and applied
this technique in a number of different situations. The power of this approach
is apparent from \cite{Obj2,Has}.

A Kodaira (complex) surface is a compact quotient of $N_2$ endowed with an
integrable left-invariant complex structure, and actually admits a holomorphic
symplectic structure. It corresponds to $g=1$ in

\proc{Theorem}{FGM91} A compact real surface $U$ of genus $g\ge1$ with a
nowhere zero 2-form and symplectomorphism $\varphi$ fixing a non-zero class in
$H^1(U,\Z)$ defines a circle bundle $E\to(U\t[0,1])/\varphi$ that is symplectic
and generally non-K\"ahler.

\n A further generalization of this construction is shown to account for all
symplectic manifolds with a free $S^1$ action.\mb

Some tools are available to compute the Dolbeault cohomology of an invariant
complex structure on $N/\Ga$ \cite{Obj1,CF,PP}, but not in the most general
situation. It is therefore natural to study the convergence of the Fr\"olicher
spectral sequence $E^{p,q}_r$, that relates the Dolbeault and de\,Rham
groups. Whereas any K\"ahler manifold (and complex surface) satisfies
$E^{p,q}_2 =E^{p,q}_\infty$, few explicit examples of higher order degeneration
were known prior to

\proc{Proposition}{CFG91} There exist invariant complex structures on
nilmanifolds of real dimension $d$ such that (i) $d=8$ and $E_2\ne E_\infty$,
and (ii) $d=12$ and $E_3\ne E_\infty$.

\n Semisimple examples, but in somewhat higher dimensions, were provided by
\cite{Pit}.\mb

Gray himself introduced me to the classification problem for invariant complex
structures on nilmanifolds, and I tackled the 6-dimensional case \cite{S6}
building on results of \cite{CFGU}. Invariant symplectic structures are
classified in \cite{GK}, again from a Lie algebra perspective. Combining the
various approaches yields

\pro{Theorem} There are 34 isomorphism classes of real 6-dimensional nilpotent
Lie algebras of which\y\hspace{10pt}\begin{tabular}{l} 
15 admit both complex and symplectic structures,\\
3 admit complex but not symplectic structures,\\
11 admit symplectic but not complex structures,\\
5 admit neither complex nor symplectic structures.\end{tabular}

These examples provide a forum for the investigation of the Gray-Hervella
classes discussed in \S2. An inner product on a 6-dimensional nilpotent Lie
algebra gives rise to a metric $g$ on an associated nilmanifold $M=N/\Gamma$.
The corresponding set of almost-Hermitian structures on $M$ is isomorphic to
\be\frac{SO(6)}{U(3)}\cong\CP^3,\label{CP3}\ee a point of which defines an
almost complex structure $J$. For each $S\subseteq\{1,2,3,4\}$ there is a
corresponding subset \[Z_S=\{J\in\CP^3:\na J\in\bigoplus\limits_{i\in S}
\cW_i\}.\] Recall that the component of $\na J$ in $\cW_4$ is affected by a
conformal change in the metric. In the lecture (available on my homepage), I
coloured points of a tetrahedron representing \rf{CP3} red, green, blue to
measure the respective components in $\cW_1,\cW_2,\cW_3$, and displayed the
resulting spectra for the example below. Whilst this was partly light-hearted,
it emphasized that the resulting lattice of 16 subsets of \rf{CP3} is subject
to non-trivial constraints arising from the non-existence of a K\"ahler metric
on $M$.\mb

\ex. A full description of the classses for a standard metric on the Iwasawa
manifold $H/\Ga$ has the following key features. The complex Heisenberg group
$H$ possesses (almost by definition) a bi-invariant complex structure $J_0$
that defines an `origin' in $\CP^3$. The torsion $\na J$ has zero
$\cW_4$-component if and only if $J$ lies in the union of distinct hyperplanes
$F,F'$ in $\CP^3$, so $Z_{\{1,2,3\}}=F\cup F'$. Let $L=F\cap F'\cong\CP^1$. 
Then

\vfil\eject

\n(i) $F\cup F'$ contains all 15 proper classes, and $J_0\in F\setminus F'$.

\n(ii) The set $Z_{\{3,4\}}=Z_{\{3\}}$ of Hermitian structures equals
$\{J_0\}\sqcup L$.

\n(iii) The set $Z_{\{2\}}$ of `almost K\"ahler' or compatible symplectic
structures is a 3-sphere in $F$ separating $J_0$ and $L$.

\n More details are contained in \cite{AGS01}. Further examples of
nilmanifolds of dimension 6 with $b_1=4$ exhibit isolated Hermitian
structures (see the remark after \rf{Wpm}), and illustrate the way the
homotopy type of a subset $Z_S$ can depend on the choice of metric.\mb

\sss{Manifolds whose holonomy is not a subgroup of $Sp(2)Sp(1)$}

The note \cite{Gnote} helped to enhance interest in a subject that had already
been founded by Berger, Bonan, Ishihara and others. It contains the important
result that a quaternionic submanifold of a quaternion-K\"ahler (or
hyperk\"ahler) manifold is totally geodesic. Whilst this put paid to hopes of
developing a theory of quaternionic submanifolds, it took many years to realize
that the most appropriate way of generating new quaternionic manifolds is
instead a quotient construction \cite{HKLR}.

A quaternion-K\"ahler manifold of dimension $4n$ can be characterized by the
existence of a parallel 4-form with stabilizer $Sp(n)Sp(1)$. It is known that
closure of the 4-form $\W$ implies covariant constancy of $\W$ in dimension
$4n\ge12$ \cite{Swa}, but this has left an intriguing open question in
dimension 8. Below, we establish the existence of a compact 8-manifold with a
closed non-parallel 4-form with stabilizer $Sp(2)Sp(1)$. This is a quaternionic
analogue of the symplectic manifolds mentioned in the previous section that do
not admit a K\"ahler metric.

Motivation comes from the (possibly) more familiar theory of $\Spin7$
holonomy. Consider the 2-forms \be\ba{rcl} \w_1 &=&
e^{13}+e^{57}+e^{24}+e^{68},\\\w_2 &=& e^{15}-e^{37}+e^{26}-e^{48},\\\w_3 &=&
e^{17}+e^{35}+e^{28}+e^{46}.\ea\label{www}\ee on $\R^8$, whose quaternionic
structure is defined by first dividing the coordinates into `odd' and
`even'. If we set \be\W_\pm=\w_1\we\w_1+\w_2\we \w_2
\pm\w_3\we\w_3,\label{pm}\ee then $\W_+$ has stabilizer $Sp(2)Sp(1)$, and
$\W_-$ has stabilizer $\Spin7$ \cite{BH}.

It is instructive to carry out `dimensional reduction' by regarding the 1-forms
$e^7,e^8$ as constant. Let \be\si=-e^{12}+e^{34}+e^{56}.\label{si}\ee Then
\[\ft12\W_-=\ft12\si^2+\phi\we e^7+\psi\we e^8-\si\we e^{78},\]
where \[\phi+i\psi=(e^1+ie^2)\we(e^3+ie^4)\we(e^5+ie^6).\] The equation
$d\W_-=0$ is therefore satisfied on $\R^6\times\R^2$ or a compact manifold of
the form $M^6\t T^2$ if $M^6$ has a symplectic form $\si$ and certain
compatible closed 3-forms $\phi,\psi$. The stabilizer of both $\phi$ and $\psi$
is the same subgroup $SL(3,\C)$ of $GL(6,\R)$; such structures have been
investigated in \cite{Hit}. Since $\si\we\phi=0=\si\we\psi$, the 3-form
$\phi+i\psi$ determines a complex structure on $M$ for which $\si$ has type
$(1,1)$. It follows that $M$ is necessarily K\"ahler and so Calabi-Yau (or
flat). In this case, $M\times T^2$ has holonomy in
$SU(3)\times\{e\}\subset\Spin7$.

The `plus' case in \rf{pm} yields more flexibility:

\pro{Theorem} There exists a closed 4-form with stabilizer $Sp(2)Sp(1)$ on a
compact nilmanifold of the form $M^6\times T^2$. The associated Riemannian
metric $g$ is reducible and is not therefore quaternion-K\"ahler.

\n\textit{Proof.} By analogy to \rf{si}, let \[\tau=e^{12}+e^{34}+e^{56}.\] Then
\[\ft12\W_+=-\ft12\tau^2+\a\we e^7+\b\we e^8-\tau\we e^{78},\] where \be
\ba{rcl}\a&=&3e^{135}+e^{146}+e^{236}+e^{245},\\\b&=&3e^{246}+e^{235}+e^{136}+
e^{145}.\ea\label{ab}\ee As before, we seek a 6-manifold $M$ with a symplectic
form $\tau$ and closed 3-forms $\a,\b$, which satisfy
$\tau\we\a=0=\tau\we\b$. This time, the stabilizer of each of $\a,\b$ is
isomorphic to $SL(3,\R)\times SL(3,\R)$ but (because of the factor 3) these
stabilizers do not coincide. Indeed, the overall structure group
$Sp(2)Sp(1)\cap SO(6)$ of $M$ is $SO(3)$ acting diagonally on $\R^3\op\R^3$.

Take $M$ to be a 6-dimensional nilmanifold associated to the Lie algebra $\g=
\al e^2\ar\op\al e^1,e^3,e^4,e^5,e^6\ar$ determined by the relations \be
\left\{\ba{l} de^i=0,\q i=1,2,3,5,\\de^4=e^{15},\\de^6=e^{13}.\ea\right. 
\label{M6}\ee The fact that
the structure constants are rational guarantees the existence of a lattice
$\Gamma$ in the associated Lie group $G$ for which $M=G/\Gamma$ is compact.

Equations \rf{M6} imply that $d\tau=0$. Moreover, of the 8 simple 3-forms in
\rf{ab}, all are closed except $e^{246}$. Let $\hW$ be the 4-form obtained from
$\W$ by substituting \be\ba{rcl} e^1 \lar e^1+\sqrt3 e^2,\\e^3 \lar e^3-\sqrt3
e^4.\ea\label{sub}\ee The stabilizer of $\hW$ is still isomorphic to
$Sp(2)Sp(1)$, though the 1-form $e^2$ will no longer be covariant constant
relative to the new metric, leaving only $\langle e^7,e^8\rangle$ as a direct
summand. Observe that \rf{sub} leaves $\tau$ unchanged. Its effect on (1) is
given by \[\ba{rcl} \a\lar\a - 9e^{245}-3\sqrt3e^{145}+3\sqrt3 e^{235},\\
\b\lar\b - 3e^{246}-\sqrt3e^{146}+\sqrt3 e^{236}.\ea\] The offending term
$3e^{246}$ has been eliminated from $\a$ at the expense of adding only closed
3-forms.  It follows that $\hW$ is closed.\hfill QED\vs

Note that the tensor $\hW$ is irrational relative to a basis for which the
structure constants of $\g$ are rational. The fact that $\hW$ is not parallel
also follows from verifying that the ideal generated by \rf{www} is not closed
under exterior differentiation. A similar reduction to 7 (rather than 6)
dimensions will allow such structures to be built up from 7-manifolds with a
certain type of $G_2$ structure arising from the subgroup $Sp(2)Sp(1)\cap
SO(7)$ isomorphic to $SO(4)$.\vs\vs

\n\textbf{Acknowledgments.} Untold thanks are due to the Organizing Committee
who worked so hard to make the conference in Bilbao a success. I am also
grateful to Mary Gray and Wolfram Research for their generosity in providing me
with versions of Mathematica to use in my lecture and beyond.

\bibliographystyle{apalike}

{\small 
\bibliography{aps}
}

\enddocument